\documentclass[12pt]{article}
\usepackage{amsmath}
\usepackage[noend]{algpseudocode}
\usepackage{algorithm}
\usepackage{amssymb}
\usepackage{graphicx}
\usepackage{hyperref}
\usepackage{float}

\newcommand{\reals}{{\mbox{\bf R}}}

\newcommand{\symm}{{\mbox{\bf S}}}  
\newcommand{\lambdamax}{{\lambda_{\rm max}}}
\newcommand{\ie}{{\it i.e.}}
\newcommand{\eg}{{\it e.g.}}
\newcommand{\argmin}{\mathop{\rm argmin}}

\newcommand{\prox}{\mathop{\bf prox}}
\newcommand{\dom}{\mathop{\bf dom}}
\newcommand{\diag}{\mathop{\bf diag}}
\newcommand{\BIT}{\begin{itemize}}
\newcommand{\EIT}{\end{itemize}}
\newcommand{\BEQ}{\begin{equation}}
\newcommand{\EEQ}{\end{equation}}
\newcommand{\BEAS}{\begin{eqnarray*}}
\newcommand{\EEAS}{\end{eqnarray*}}

\newcounter{algorithmctr}[section]
\renewcommand{\thealgorithmctr}{\thesection.\arabic{algorithmctr}}
\newenvironment{algdesc}%
{\refstepcounter{algorithmctr}
\begin{list}{}{%
   \setlength{\rightmargin}{0\linewidth}%
   \setlength{\leftmargin}{.05\linewidth}}%
   \rmfamily\small
   \item[]{\setlength{\parskip}{0ex}\hrulefill\par%
    \nopagebreak{\bfseries\textsf{Algorithm \thealgorithmctr~}}}}%
{{\setlength{\parskip}{-1ex}\nopagebreak\par\hrulefill}
\end{list}}

{\refstepcounter{algorithmctr}
\begin{list}{}{%
   \setlength{\rightmargin}{0\linewidth}%
   \setlength{\leftmargin}{.05\linewidth}}%
   \rmfamily\small
   \item[]{\setlength{\parskip}{0ex}\hrulefill\par%
    \nopagebreak{\bfseries\textsf{Algorithm}}}}%
{{\setlength{\parskip}{-1ex}\nopagebreak\par\hrulefill}
\end{list}}

\begin{document}

\title{Polyak Minorant Method\\ for Convex Optimization}
\author{Nikhil Devanathan \and Stephen Boyd}
\date{\today}
\maketitle

\begin{abstract}
In 1963 Boris Polyak suggested a particular step size for gradient
descent methods, now known as the Polyak step size, 
that he later adapted to subgradient methods.
The Polyak step size requires knowledge
of the optimal value of the minimization problem,
which is a strong assumption but one that holds for
several important problems.
In this paper we extend Polyak's method to handle constraints 
and, as a generalization of subgradients, general
minorants, which are convex functions that
tightly lower bound the objective and constraint functions.
We refer to this algorithm as the Polyak Minorant Method (PMM).
It is closely related to cutting-plane and bundle methods.
\end{abstract}

\clearpage

\section{Introduction}

\subsection{The problem}
We consider the convex optimization problem
\BEQ\label{e-prob}
\begin{array}{ll}
\mbox{minimize}   & f_0(x) \\
\mbox{subject to} & f_i(x) \leq 0, \quad i=1, \ldots, m\\
& Ax=b,
\end{array}
\EEQ
with variable $x\in \reals^n$,
where $f_0: \Omega \to \reals$ and $f_i: \reals^n \to \reals$, $i=1, \ldots, m$
are closed proper convex functions, $A\in \reals^{p \times n}$, and $b\in \reals^p$.
We can have $m=0$ (no inequality constraints) or $p=0$ (no equality constraints).

We let $\mathcal F$ denote the set of feasible points for \eqref{e-prob}.
For $\rho>0$, $\mathcal F_\rho$ will denote the $\rho$-violated
constraint set,
\[
\mathcal F_\rho = \{ x \mid f_i(x)\leq \rho,~i=1, \ldots, m,~Ax=b\}.
\]
We will assume that $\mathcal F_\rho \subseteq \Omega$ for some
$\rho>0$, which means that when $x$ is $\rho$-close to feasible, 
the objective $f_0(x)$ is defined.
We assume that the optimal value
\[
f^\star = \inf \{f_0(x) \mid f_i(x)\leq 0,~i=1,\ldots, m,~Ax=b\}
\]
is finite and achieved,
\ie, there exists at least one optimal point $x^\star$, with $f_i(x^\star) \leq 0$
for $i=1,\ldots,m$ and $f_0(x^\star)=f^\star$.
We define $\mathcal{X}^\star$ to be the set of optimal points.
For our convergence proofs, we will also assume that
$f_i$ for $i=0,\ldots, m$ are Lipschitz continuous with constant $G$.

We define the (maximum) violation
at a point $x\in \Omega$ that satisfies $Ax=b$ as
\BEQ\label{e-tv}
v(x) = \max\{f_0(x)-f^\star, f_1(x),\ldots,f_m(x)\},
\EEQ
and take $v(x)=\infty$ when $x \not\in \Omega$ or $Ax \neq b$.
The violation is zero if and only if $x$ is a solution of \eqref{e-prob}.
An algorithm solves the problem \eqref{e-prob} if it produces a sequence $x^k$
with $v(x^k) \to 0$.

\subsection{Known optimal value} \label{s-problems}
Like Polyak's original method, 
the algorithm we present in this paper assumes knowledge of $f^\star$.
Although requiring that $f^\star$ is known before solving the problem is restrictive,
there are common generic cases where it holds.

\paragraph{Feasibility problems.}
A feasibility problem is the special case of \eqref{e-prob} with $f_0=0$,
\BEQ\label{e-f-prob}
\begin{array}{ll}
	\mbox{minimize}   & 0 \\
	\mbox{subject to} & f_i(x) \leq 0, \quad i=1, \ldots, m\\
	& Ax=b,
\end{array}
\EEQ
with variable $x\in \reals^n$, and
$f_i$, $A$, and $b$ as in \eqref{e-prob}.

\paragraph{Primal-dual problems.}
A primal-dual problem includes both primal and dual variables and 
constraints, and includes the duality gap (the difference
of the primal objective and the dual objective) as the objective
or as a constraint (that it is zero).
Such a problem has known objective value $0$ when strong duality holds and 
the primal problem has a solution.
As a specific example consider the primal and dual cone programs
\BEQ\label{e-conepair}
\begin{array}{ccc}
\begin{array}{ll}
\mbox{minimize}   & c^T u \\
\mbox{subject to} & Au=b\\
& u\in \mathcal{K},
\end{array}
& \qquad &
\begin{array}{ll}
\mbox{maximize}   & b^T v \\
\mbox{subject to} & c - A^Tv = s\\
& s\in \mathcal{K}^*\
\end{array}
\end{array}
\EEQ
with variables $u\in \reals^n$, $s\in \reals^n$, $v\in \reals^p$,
where $\mathcal K \subseteq \reals^n$
is a closed convex cone and $\mathcal{K}^*$ is its dual cone,
and $c\in \reals^{n}$, $A\in \reals^{p\times n}$, and 
$b\in \reals^{p}$ are parameters.
Expressing the condition that the duality gap $c^Tu-b^Tv$ is zero as 
as a linear constraint, we arrive at the primal-dual feasibility problem
\BEQ\label{e-pdcone-pmm}
\begin{array}{ll}
\mbox{minimize}   & 0 \\
\mbox{subject to} & d_{\mathcal{K}}(u)\leq 0, \quad d_{\mathcal{K}^*}(s)\leq 0\\
& \left[\begin{array}{c}
	s \\
	0 \\
	0 \\
\end{array}\right]=
\left[\begin{array}{ccc}
	0 & -A^T & c\\
	A & 0 & -b\\\
	-c^T & b^T & 0
\end{array}\right]\left[\begin{array}{c}
	u \\
	v \\
	1
\end{array}\right],
\end{array}
\EEQ
where $d_{\mathcal{K}}$ is the $\ell_2$ distance to $\mathcal{K}$ and 
$d_{\mathcal{K}^*}$ is the distance to $\mathcal{K}^*$. 
This has the form \eqref{e-prob} with variable $x=(u,v,s)$ 
and known optimal value $f^\star = 0$.

\subsection{Pointwise lower bounds and minorants}
\paragraph{Pointwise lower bound.}
Suppose $f: \Omega \to \reals$, with $\Omega \subseteq \reals^n$.
We say that $\hat f: \reals^n \to \reals \cup \{\infty\}$ is a pointwise 
lower bound (PLB) on $f$ if $\hat f(z) \leq f(z)$ for all $z\in \Omega$. 
We write this as $\hat f \leq f$.

\paragraph{Minorant.}
We will use the basic idea of a minorant of a convex function at a point.
Suppose $f:\Omega \to \reals$ is a closed convex function.
We say that $\hat f : \reals^n \to \reals \cup \{\infty\}$ is a minorant
of $f$ at $z\in \Omega$ if the following hold:
\BIT
\item $\hat f$ is closed convex;
\item $\hat f \leq f$, \ie, $\hat f$ is a PLB on $f$;
\item $\hat f(z) = f(z)$, \ie, the lower bound is tight at the point $z$.
\EIT
We write $\hat f(x)$ as $\hat f(x;z)$ to indicate that $\hat f$ is a minorant
at the point $z$.

The simplest minorant of $f$ at $z$ is affine,
\BEQ\label{e-aff-minorant}
\hat f(x;z)=f(z)+g^T(x-z), \qquad
g\in \partial f(z),
\EEQ
with $\partial f(z)$ denoting the subdifferential of $f$ at $z$.
At the other extreme, we can take $f$ as its own minorant, $\hat f(x;z) = f(x)$.
All other minorants are in between these two, in the sense that they are
pointwise between $f$ and at least one affine minorant defined by a subgradient.
We will say more about minorants and how to construct them
in \S\ref{s-minorants}.

If $f$ has Lipschitz constant $G$,
we will also assume that any minorant inherits the same Lipschitz constant.
(Some minorants do not inherit the Lipschitz constant,
but all the typical methods for constructing minorants do have this
property.)

\paragraph{Other notation.}
Many authors have used different names for what we call a PLB and a minorant.
Some papers define (what we call) a minorant to be (what we call) a PLB
\cite{parshakova2022implementation,drusvyatskiy2020convex}.
Other authors use minorant to be more restrictive, for example
affine \cite{lin2018level,mclinden1978affine}.
In \cite{jongbloed1998the}, a
minorant is defined as the largest convex function whose epigraph contains
a set of points.
Throughout this paper, we use the definitions of PLB and minorant above.

\subsection{This paper}
In this paper, we develop a method that solves \eqref{e-prob} where the only access
to the objective and constraint functions is via minorants.  That is,
we can find a minorant of $f_i$, $i=1,\ldots, m$ 
at any $z$, and of $f_0$ for any $z\in \Omega$.
Our method is inspired
by and is an extension of the subgradient method with Polyak step size. As we will
discuss below, it is closely connected to many other methods including subgradient
methods, cutting-plane methods, and bundle methods. Much like bundle and cutting-plane
methods, each iteration of our method requires solving a relatively simple
convex optimization problem.  One benefit of the method is that
it has no parameters that need to be tuned.
To honor Boris Polyak, we name the method the Polyak minorant method (PMM).

\subsection{Prior work}
\paragraph{Subgradient methods.}
Boris Polyak introduced the Polyak gradient step for minimizing continuously
differentiable functionals in \cite{polyak1963gradient}. Following the development
of subgradient methods for non-differentiable optimization
\cite{shor2012minimization, rockafellar1981TheTO}, Polyak adapted his gradient
step to the subgradient case \cite{polyak1987intro}. Since then various
extensions of the Polyak subgradient method have been developed, including
Polyak step variations for stochastic gradient descent 
\cite{loizou2021stochastic,abdukhakimov2023stochastic,berrada2021comment,
gower2022cutting,li2023sp2,prazeres2021stochastic,orvieto2022dynamics},
Polyak-like steps that do not require knowledge of $f^\star$ 
\cite{hazan2022revisiting,you2022two}, a Polyak-like step for
momentum-accelerated gradient descent 
\cite{wang2023generalized,barre2022complexity,goujaud2022quadratic}, a
Polyak step method for mirror descent \cite{you2022minimizing}, a Polyak-like
step for convex problems with box constraints \cite{cheng2012active},
and a Polyak step for weakly convex functions \cite{davis2018subgradient}.
We will see that PMM reduces to the subgradient method with Polyak step sizes
when there are no constraints and a subgradient-based affine minorant is used.

\paragraph{Cutting-plane methods.}
Cutting-plane methods originate from the works \cite{cheney1959newton} 
and \cite{kelley1960cutting}. These methods solve convex
problems by iteratively shrinking a polygonal superset of the optimal set.
This shrinking is done using a cutting plane in each iteration, \ie,
a halfspace known to contain the optimal set \cite{boyd2007localization}.   
Cutting-plane methods differ in how they choose the next iterate; see, \eg,
the survey \cite{sra2012optimization}.
Cutting-plane methods can be viewed as iteratively refining a 
piecewise-affine minorant on the objective and constraints.  
Conversely, PMM can be thought of as a 
cutting-plane method, when the minorants are piecewise affine.

\paragraph{Bundle methods.}
Bundle methods extend cutting-plane methods in two ways.
First, the next iterate is found by minimizing a minorant plus 
an additional (typically quadratic) stabilization term.
Second, bundle methods include logic that only updates the 
current iterate if a sufficient descent condition holds.
The original and most common bundle method
is the proximal bundle method \cite{kiwiel1990proximity, frangioni2020standard}.
Alternatively,
the level bundle method \cite{lemarechal1995new, frangioni2020standard}
projects the current point onto the sublevel set of a minorant,
exactly as PMM does.
Yet another variation is
the trust-region bundle method \cite{marsten1975boxstep, frangioni2020standard}.
A history of bundle methods can be found in
\cite[Ch.~XIV, XV]{hiriart1996convex}. 
We refer to \cite{parshakova2022implementation}
for a more thorough review of modern bundle literature.
PMM is structurally similar to a level-set bundle method, but lacks a sufficient 
descent condition, and admits any minorant instead of only cutting-plane minorants.

\section{Polyak minorant method}
We let $x^k \in \reals^n$ denote the $k$th iterate of PMM
for $k=1,2, \ldots{}$.
PMM maintains a PLB on the objective and constraint functions,
which vary with iterations, denoted $\hat f_i^k$, $i=0,\ldots, m$.
The constraint function lower bounds $\hat f_i^k$, $i=1, \ldots, m$,
are always minorants of $f_i$ at $x^k$,
while the PLB $\hat f_0^k$ is a minorant only for some iterations.
We define
\BEQ\label{e-hat-sublevel-set}
\mathcal X^k = \{ x \mid \hat f_0^k(x) \leq f^\star,~
\hat f_i^k(x) \leq 0,~i=1,\ldots, m, ~Ax=b \}.
\EEQ
Since $\hat f_i^k \leq f_i$
for $i=1, \ldots, m$ and $\hat f_0^k \leq f_0$, we have that
$\hat f_i^k(x^\star) \leq f_i(x^\star)\leq 0$ for $i=1, \ldots, m$
and $\hat f_0^k(x^\star) \leq f_0(x^\star) = f^\star$, so $x^\star \in \mathcal X^k$.

The next iterate $x^{k+1}$ is the projection of $x^k$ onto $\mathcal X_k$,
\ie,
\[
x^{k+1} = \Pi_{\mathcal X^k} (x^k) =
\argmin_{x \in \mathcal X^k} \|x-x^k\|_2.
\]
PMM is summarized in algorithm \ref{a-pmm}.

\begin{algdesc}{\sc Polyak minorant method} \label{a-pmm}
{\footnotesize
\begin{tabbing}
{\bf given} $x^1 \in \reals^n$, optimal value $f^\star$, tolerance $\epsilon >0$.
\\*[\smallskipamount]
{\bf for $k=1,2,\ldots$} \\
\qquad \= 1.\ \emph{Constraint minorants.}
Find minorants $\hat f_i^k$ of $f_i$ at $x^k$ for $i=1, \ldots, m$.\\
\> 2.\ \emph{Objective minorant.}\\
\>\qquad\= If $x^k \in \Omega$, find a minorant $\hat f_0^k$ of $f_0$ at $x^k$.\\
\>\> Else, set $\hat f_0^k$ to be any PLB on $f_0$.\\
\> 3.\ \emph{Update.} $x^{k+1} = \Pi_{\mathcal X^k}(x^k)$.\\
\> 4.\ \emph{Check stopping criterion.} Stop if $v(x^{k+1}) \leq \epsilon$.
\end{tabbing}}
\end{algdesc}

\paragraph{Comments.}
Recall that $\dom f_i = \reals^n$ for $i=1, \ldots, m$,
so we can always find constraint minorants in step~1.
In step~2, if $x^k \in \Omega =\dom f_0$, $\hat f_0^k$ is a minorant of $f_0$ at $x^{k}$.
If $x^k \not\in \Omega$, $\hat f_0^k$ is any PLB, such as
the constant function $\hat f_0^k=f^\star$, or a minorant of $f_0$
found in any previous iteration.

The PMM method is generic since we have not specified
what minorants to use. We will discuss many different types of
minorants in \S\ref{s-minorants}.  Depending on the minorants used,
PMM can be considered a subgradient-type method, a cutting-plane method, or a
level-set bundle method \cite{van2016inexact}.

\paragraph{Connection to proximal operator.}
Define
\[
\hat F^k(x) = \hat f_0(x) + \mathcal I(\hat f_i(x) \leq 0,~i=1, \ldots, m,~Ax=b),
\]
where $\mathcal I$ is the $\{0,\infty\}$-indicator function.
The projection of $x^k$ onto $\mathcal X^k$ also minimizes
$t^k\hat F^k(x) + (1/2) \|x-x^k\|_2^2$ for some $t^k>0$.
In other words, we have $x^{k+1}=\prox_{t^k\hat F^k}(x^k)$, where
$\prox$ is the proximal operator \cite{parikh2014proximal}.
In standard proximal operator methods, $t^k$ is specified.
In our case, however, $x^{k+1}$ is given as a projection, and we
determine $t^k$ only after the update is computed.

\paragraph{Polyak subgradient method.}
Consider the special case $\Omega=\reals^n$ and $m=p=0$.
We use the subgradient-based affine minorant \eqref{e-aff-minorant} for $f_0$.
The set $\mathcal X^k$ is the halfspace
$\{x \mid f_0(x^k) + (g^k)^T(x-x^k) \leq f^\star \}$,
where $g^k \in \partial f_0(x^k)$.
The projection in step~3 is then
\[
x^{k+1} = x^k - \frac{f_0(x^k)-f^\star}{\|g^k\|_2^2} g^k,
\]
which coincides with the subgradient method with Polyak's step size.
(This assumes $g^k \neq 0$; if $g^k=0$, we can terminate since $x^k$ is optimal.)
So PMM generalizes the subgradient method with Polyak step size.

\paragraph{Alternating-update PMM.}
We mention one
simple variation of PMM in which the projection is replaced with
a projection onto the objective sublevel set or a projection onto
the constraint minorant sublevel set. We define
\BEAS
\mathcal{X}_0^k&=&\{x\mid \hat f_0^k(x)\leq f^\star\},\\
\mathcal{X}_1^k&=&\{x\mid \hat f_i^k(x) \leq 0, ~i=1,\ldots, m, ~\ Ax=b\}.
\EEAS
We replace the projection in step~3 of PMM \ref{a-pmm} with projection onto
$\mathcal{X}_0^k$ for even $k$ and $\mathcal{X}_1^k$ for odd $k$. 
This alternating-update variation is very similar to existing
alternating-update subgradient methods for constrained convex optimization
\cite{polyak1967general, beck2010comirror, lan2020algorithms}. This 
modified PMM converges under the same assumptions as the original PMM;
we provide a brief proof of convergence in \S\ref{s-alt-convergence}.

\subsection{Convergence proof} \label{s-convergence}
Here we give a short convergence proof for PMM, \ie, we 
show that $v(x^k)\to 0$ as $k\to \infty$, which implies that
the stopping criterion is eventually satisfied.
(A very similar proof can be constructed to show that the alternating-update 
version also converges.)
We give the proof not because it is novel, but because it is short
and simple.
It uses basic convex analysis and standard ideas that trace back to the subgradient 
methods of the 1960s.

Since $x^{k+1}$ is the projection of $x^k$ onto $\mathcal X^k$ 
(which contains $x^\star$), we have
\[
(x^k-x^{k+1})^T(x^\star-x^{k+1})\leq 0.
\]
It follows that
\BEAS
\|x^{k+1} - x^\star\|_2^2 &= & \|x^k-x^\star\|_2^2-\|x^k-x^{k+1}\|_2^2
+2(x^k-x^{k+1})^T(x^\star-x^{k+1})\\
&\leq & \|x^k-x^\star\|_2^2-\|x^k-x^{k+1}\|_2^2.
\EEAS
This shows that the algorithm is Fej\'er monotone, \ie, each iteration does
not increase the distance to any optimal point.
Iterating the inequality above yields
\[
\sum_{k=1}^\infty \|x^k-x^{k+1}\|_2^2 \leq \|x^1-x^\star\|_2^2,
\]
which shows that
\BEQ\label{e-norm-diff}
\|x^k-x^{k+1}\|_2\to 0.
\EEQ
We use the Lipschitz continuity of $\hat f_i^k$,
the fact that $x^{k+1} \in \mathcal X^k$
(which implies $\hat f_i^k(x^{k+1}) \leq 0$), and that $\hat f_i^k$ is a minorant
of $f_i$ at $x^k$ to conclude that
\BEQ\label{e-norm-diff-f1}
G\|x^{k+1}-x^k\|_2\geq \hat f_i^k(x^k) - \hat f^k_i(x^{k+1}) \geq
f_i(x^k)
\EEQ
for $i=1,\ldots,m$.

Combining \eqref{e-norm-diff} and \eqref{e-norm-diff-f1}, we see that
$\max\{ f_i(x^k) , 0\} \to 0$ as $k \to \infty$, for $i=1,\ldots, m$.
In other words, the iterates are eventually almost feasible.
It follows that there exists a $K$ for which
$f_i(x^k)\leq \rho$, $i=1,\ldots, m$, for all $k \geq K$,
which implies $x^k \in \mathcal F_\rho$, and thus $x^k \in \Omega$.
This in turn implies that $\hat f_0^k$ is a minorant of $f_0$ at $x^k$
for $k \geq K$.

We now show that $f_0(x^k) \to f^\star$.
For $k \geq K$, a similar argument as above for
$f_i$, $i=1,\ldots, m$, gives
\BEQ\label{e-norm-diff-f0}
G\|x^{k+1}-x^k\|_2\geq \hat f_0^k(x^k) - \hat f^k_0(x^{k+1}) \geq
\hat f_0(x^k) - f^\star = f_0(x^k) - f^\star.
\EEQ
(Here we use $\hat f_0(x^{k+1}) \leq f^\star$ and
$\hat f_0(x^{k}) = f_0(x^k)$.)
Combined with \eqref{e-norm-diff}, we deduce that $f_0(x^k) \to f^\star$.
It follows that $v(x^k) \to 0$, so the stopping criterion is eventually 
satisfied.

\paragraph{Comments.}
The convergence proof shows that we can relax several assumptions.
For example, the assumption of Lipschitz
continuity can be relaxed by requiring it to hold for the minorants, and only on the
bounded set defined by $\|x-x^1\|_2 \leq \|x^\star -x^1\|_2$.
(We do not know the righthand side here, but we never use
the Lipschitz constant in the algorithm.)

While we have assumed above that all constraint functions have domain
$\reals^n$, our proof shows that it suffices for just one to be defined on
all $\reals^n$, with the other constraint functions playing a similar role to
the objective, \ie, they are defined only when the one constraint is nearly
satisfied.  Instead of a minorant, we take $\hat f_i$ to be any PLB
for $f_i$ for the constraints with $x^k \not\in \dom f_i$, as we do
above for the objective.

\subsection{Convergence with alternating update}\label{s-alt-convergence}
We follow the steps of the proof in \S\ref{s-convergence} to 
show that $f_0(x^k) \to f^\star$ as $k\to \infty$ for the alternating-update 
variant of the PMM. By construction, $Ax^k=b$. Since $x^\star\in X_0^k$ and 
$x^\star\in X_1^k$, we deduce that for both odd and even $k$,
\[
(x^k-x^{k+1})^T(x^\star-x^{k+1})\leq 0.
\]
It follows that 
\[
\|x^{k+1} - x^\star\|_2^2 \leq\|x^k-x^\star\|_2^2-\|x^k-x^{k+1}\|_2^2.
\]
Iterating these inequalities yields that
\[
\sum_{k=1}^\infty \|x^k-x^{k+1}\|_2^2 \leq \|x^1-x^\star\|_2^2.
\]
It now follows that $\|x^k-x^{k+1/2}\|_2\to 0$. We conclude using
\eqref{e-norm-diff-f1}, \eqref{e-norm-diff-f0}, and the Lipschitz 
continuity of $\hat f_i$ for $i=0,\ldots m$ that 
$f_0(x^k) \to f^\star$ as $k\to \infty$.

\subsection{Convergence with sharpness}
A convex function $\phi$ on $\Omega$, with $\phi^\star = \inf_{x \in \Omega}\phi(x)$
finite,
is $\mu$-sharp if for some $\mu> 0$ and all $x\in\Omega$
\[
\phi(x)-\phi^\star\geq \mu \inf_{x^\star\in \mathcal{X}^\star}\|x-x^\star\|_2.
\]
It is a well-known result (see \cite{polyak1969minimization,shor1973convergence,
goffin1977rates,davis2018subgradient}) that for minimizing a $G$-Lipschitz 
$\mu$-sharp closed proper convex function, the original Polyak subgradient 
step produces iterates that converge linearly:
\[
\inf_{x^\star\in \{x\mid \phi(x)\leq \phi^\star\}}\|x^{k+1}-x^\star\|_2\leq 
\sqrt{1-\frac{\mu^2}{G^2}}\inf_{x^\star\in \{x\mid f(x)\leq f^\star\}}
\|x^{k}-x^\star\|_2.
\]

We can derive a similar convergence result for PMM.
Define 
\[
h(x)=\max(f(x)-f^\star,f_1(x),\ldots,f_m(x))
\]
on $\Omega$. The function $h$ is $G$-Lipschitz, closed, proper,
and convex, and its set of minimizers is $\mathcal{X}^\star$.
Now suppose that $h$ is $\mu$-sharp.
From \eqref{e-norm-diff-f1}, \eqref{e-norm-diff-f0}, and 
sharpness, we have
\[
G\|x^{k+1}-x^{k}\|_2\geq h(x^k)\geq 
\mu\inf_{x^\star\in \mathcal{X}^\star}\|x^{k}-x^\star\|_2.
\]
Letting $\Pi_{\mathcal{X}^\star}x^k$ denote the projection of $x^k$ onto
$\mathcal{X}^\star$, we have
\BEAS
\inf_{x^\star\in \mathcal{X}^\star}\|x^{k+1}-x^\star\|_2^2&\leq& 
\|x^{k+1} - \Pi_{\mathcal{X}^\star}x^k\|_2^2 \\
&\leq& \|x^k-\Pi_{\mathcal{X}^\star}x^k\|_2^2-\frac{\mu^2}{G^2}
\inf_{x^\star\in \mathcal{X}^\star}\|x^{k}-x^\star\|_2^2\\
&=&\left(1-\frac{\mu^2}{G^2}\right)
\inf_{x^\star\in \mathcal{X}^\star}\|x^{k}-x^\star\|_2^2.
\EEAS
Thus when the function $h$ is $\mu$-sharp, 
PMM attains the same linear rate of convergence
as the original Polyak subgradient method.

\subsection{Cost of an iteration}\label{s-iter}
We address here a question that could just as well be asked about cutting-plane
or bundle methods: What is the computational cost of an iteration of PMM,
specifically the projection step?
Of course, this depends very much on the minorants used.
For example, if we take the functions themselves as minorants, PMM converges
in one step, which consists of solving the problem; PMM is correct in this case,
but silly.
To be useful, carrying out the projection step should be, at a very minimum,
cheaper than solving the original problem.

Typical minorants are piecewise affine, defined as the maximum of a set of 
affine functions.  For such problems, the projection can be solved in time 
that is linear in $n$, and quadratic in the number of terms in the minorants
plus equality constraints.  In typical cases the latter number is kept 
substantially smaller than $n$ as the algorithm proceeds, using limited
memory minorants described in \S\ref{s-mem}.

While the details depend on the specific form of the minorants, we give
them here for a specific generic case, where the projection can be 
expressed as
\BEQ\label{e-qp}
\begin{array}{ll}
	\mbox{minimize}   & \|x-x^k\|_2^2 \\
	\mbox{subject to} & Fx\leq g, \quad Ax=b,
\end{array}
\EEQ
where $F\in \reals^{q \times n}$ and $g\in \reals^q$, where $q$ is the number of
terms in the piecewise affine minorants. We assume here that $q \ll n$, and show
how to solve this problem efficiently.

From the optimality condition for this quadratic program (QP), we find that 
the solution to \eqref{e-qp} has the form
\[
x^{k+1}=x^k -F^T\lambda - A^T\nu,
\]
for some dual variables $\lambda\in \reals^q$ and $\nu\in \reals^p$, with
$\lambda \geq 0$ \cite[\S 5.5.3]{boyd2004convex}. 
So we can reformulate \eqref{e-qp} using variables $\lambda$ and $\nu$ as
\BEQ\label{e-small-qp}
\begin{array}{ll}
	\mbox{minimize}   & \|F^T\lambda +A^T\nu\|_2^2 \\
	\mbox{subject to} & Fx^k-FF^T\lambda -FA^T\nu \leq g\\
	& Ax^k -AF^T\lambda -AA^T\nu=b.
\end{array}
\EEQ
This a QP with variables $(\lambda,\nu) \in \reals^{q+p}$.
We solve this (smaller) QP and then set $x^{k+1}=x^k -F^T\lambda^\star 
- A^T\nu^\star$. When $q+p \ll n$, this has far fewer variables than the 
original projection QP \ref{e-qp}.

Without exploiting any structure, the small QP \eqref{e-small-qp}
can be solved in $O((p+q)^3)$ flops
\cite[\S 11]{boyd2004convex}. 
In many cases the computational cost of the projection 
is dominated by forming the matrices
\BEQ\label{e-gram-matrices}
FF^T, \quad FA^T, \quad AA^T, 
\EEQ
which has cost $O(n(p+q)^2)$ flops.

To illustrate this, we consider a specific example 
with $n=10^6$ and $p=q=50$.
Computing the matrices \eqref{e-gram-matrices} costs $O(10^{10})$ flops, 
whereas solving the small QP \eqref{e-small-qp} costs $O(10^6)$ flops, 
which is negligible in comparison.
Carrying out this computation using CVXPY \cite{diamond2016cvxpy} and 
OSQP \cite{stellato2020osqp} on an instance of this problem 
yields results that are compatible with these rough flop counts.  
Computing the matrices requires around 0.3 seconds, 
and solving the small QP requires 0.007 seconds.
In comparison, solving the original QP directly takes around 700 seconds.
(These numbers are for a laptop with a Ryzen 9 5900HX processor.)

Similar methods to efficiently compute the projection can be used
for minorants of a more general form, \ie, not the maximum 
of $q$ affine functions.  As specific examples, the minorants
could be the sum of terms, each a maximum of affine functions,
second-order cone representable, with $q$ being something
like the total number of terms involved.
The main point here is while each iteration of PMM requires solving
a (possibly large)
convex optimization problem that does not have an analytical 
solution, this can be done efficiently provided $q$ is not too big. 

\section{Minorants} \label{s-minorants}
In this section, we look at methods for constructing minorants.
We have already mentioned some simple minorants, such as 
the affine subgradient-based minorant \eqref{e-aff-minorant}
and the function itself.
We start by mentioning simple minorants for functions that 
satisfy additional conditions, such as strongly convex or self-concordant 
\cite{nesterov1994interior}.
We then give some rules for constructing minorants, which can be 
extended to an automated method that relies on disciplined convex
programming \cite{grant2006disciplined}.

\subsection{Strongly convex and self-concordant functions}
\paragraph{Strongly convex functions.}
The subgradient-based minorant \eqref{e-aff-minorant} can be replaced
with a quadratic minorant when $f$ is 
strongly convex with parameter $\delta>0$, \ie, 
$f(x)-(\delta/2) \|x\|_2^2$ is convex.
Here the minorant is
\[
\hat f (x;z) = f(z) + g^T(x-z) + (\delta/2) \|x-z\|_2^2,
\]
where $g \in \partial f(z)$.

\paragraph{Self-concordant functions.}
As another example, suppose $f$ is
self-concordant \cite[\S 2.2.VI]{nemirovski1996interior}.  In this case,
we have the minorant
\[
\hat f(x;z) = f(z) + \nabla f(z)^T(x-z) + u - \log (1+u),
\]
where
$u = \| \nabla ^2 f(z)^{1/2} (x-z)\|_2$.
(This is convex since $u-\log (1+u)$ is increasing on $u \geq 0$.)

\subsection{Rules for constructing minorants}

\paragraph{Scaling and sum.}
If $\hat f_i(x;z)$ is a minorant for $f_i$ at $z$ and $\alpha_i \geq 0$,
then $\sum_{i=1}^M \hat f_i(x;z)$ is a minorant for 
$\sum_{i=1}^M f_i$ at $z$.
As an example, if each minorant $\hat f_i(x;z)$ is the maximum
of affine functions, then the minorant for the sum is also 
piecewise affine, with the specific form of a sum of functions,
each the maximum of affine functions.

\paragraph{Selective minorization.}
As a specific example, suppose
\[
f(x) = l(x) + \lambda r(x),
\]
where $l$ is a convex (loss) function, $r$ is convex (regularizer) function, 
and $\lambda >0$ is a parameter \cite[\S 6.3.2]{boyd2004convex}, 
\cite[\S 6.4.1]{nesterov2018lectures}.
(This function arises in regularized empirical risk minimization
problems.)
We can use the minorant
\[
\hat f(x;z) = \hat l(x;z)+ \lambda r(x),
\]
where $\hat l(x;z)$ is a minorant of $l$ at $z$.  Here we form a minorant
of the loss function, but keep the regularizer (which is its own minorant).

\paragraph{Supremum.}
Suppose $f$ is the pointwise supremum of a set of convex functions,
\[
f(x) = \sup _{\alpha \in \mathcal A} f_\alpha(x),
\]
where $f_\alpha: \Omega \to \reals$ are convex.
We will assume that the supremum is achieved for each $x$.
(When $f$ is the objective, the problem \eqref{e-prob} is then a minimax problem
\cite{boyd2004convex}.)
We can construct a minorant as follows.
Find $\alpha' \in \mathcal A$ with $f(z) = f_{\alpha'}(z)$.
Then
\[
\hat f(x;z) = f_{\alpha'}(x)
\]
is a minorant.
We can generalize this in several ways.  We can replace the righthand side
with a minorant of $f_{\alpha'}$ at $z$.
We can form a (pointwise) larger minorant as the maximum over multiple
values of $\alpha$, as long as one of them maximizes $f_\alpha(z)$.

\paragraph{Maximum eigenvalue.}
As a specific example,
we consider $f: \symm^m \to \reals$ defined as
\BEQ\label{e-smalleigenmin}
f(X) = \lambdamax(X) = \sup_{\|u \|_2=1} u^TXu,
\EEQ
where $X \in \symm^m$, the set of symmetric $m\times m$ matrices.
Using the method above, we find a minorant at $Z \in \symm^m$ 
by finding an eigenvector $v$ of $Z$, with $\|v\|_2=1$, 
associated with its maximum eigenvalue \cite{kowalewski1909power}.
This gives the minorant
\[
\hat f(X;Z) = v^TXv,
\]
which is linear. (This coincides with the minorant constructed
from the subgradient $vv^T$ of $f$ at $Z$.)
For more sophisticated (and larger) minorants, we can take
\BEQ\label{e-eigeneigenmin}
\hat f(X;Z) =\lambdamax(V^T X V),
\EEQ
or
\BEQ\label{e-diageigenmin}
\hat f(X;Z) =\max\diag(V^T X V),
\EEQ
where $V$ is an $m \times r$ matrix whose columns are orthonormal 
eigenvectors associated with the top $r$ eigenvalues of $Z$ and 
$\diag$ gives the diagonal entries of a matrix. 
These minorants reduce to \eqref{e-smalleigenmin}
when $r=1$.  For $r\geq 2$, \eqref{e-eigeneigenmin} and \eqref{e-diageigenmin}
are not equivalent, and when $r=2$, \eqref{e-eigeneigenmin} is second-order cone
representable. The minorant \eqref{e-diageigenmin} is always piecewise linear.

\subsection{DCP expressions}
Disciplined convex programming (DCP) is a system for constructing expressions
with known curvature, convex, concave, or affine
\cite{grant2006disciplined}.
Expressions are built from a library of atomic or basic functions
with known sign, monotonicity, and curvature.
These are combined in such a way that a composition rule, sufficient to
establish convexity or concavity, holds for each subexpression.
The leaves of the expression tree are constants or variables.
In addition to curvature of subexpressions, we can also track sign and
monotonicity.  Roughly speaking, we determine the signs of the zeroth, first,
and second derivatives of all subexpressions using the composition rule.

\paragraph{DCP composition rule.}
Consider the expression given by
\[
\psi_0 = \phi(\psi_1,\ldots,\psi_k),
\]
where $\phi$ is an atom and $\psi_1, \ldots, \psi_k$ are 
expressions.
The composition rule for convexity is:
$\psi_0$ is convex provided $\phi$ is convex and for each $i=1,\ldots, k$,
one of the following holds:
\begin{enumerate}
\item $\psi_i$ is affine,
\item $\psi_i$ is convex and $\phi$ is non-decreasing in argument $i$,
\item $\psi_i$ is concave and $\phi$ is non-increasing in argument $i$.
\end{enumerate}
One subtlety is that the monotonicity in conditions~2 and~3 
must be of the extended-valued function; see
\cite[\S 3.2.4]{boyd2004convex}.  There is a similar rule for concavity.

An expression is called DCP (or DCP-compliant) if every subexpression satisfies
the composition rule.
DCP is a sufficient condition for convexity or concavity.
We can think of a DCP-convex expression as a function that is syntactically 
convex, since its convexity can be established by recursively applying the
composition rule.  Note that it only requires knowing the sign, monotonicity,
and curvature of the atoms used, and not their specific values.

\paragraph{Minorant for DCP expression.}
Suppose $\psi$ is DCP-convex. 
We can construct a minorant for it by replacing every atom in it with
a minorant at its argument, provided the minorants satisfy two additional
conditions:  The minorant preserves the sign and the monotonicity of the atom.
These properties do not hold for general minorants, but they do when the 
minorants are constructed from subgradients or the other methods
described above.  Handling the sign requirement is easy; if $\phi$ is an
atom that is known to be nonnegative, we replace any minorant $\hat\phi(x;z)$
with $\max\{\hat \phi(x;z),0\}$, a minorant that is also nonnegative.

We omit the proof that this simple method produces a minorant of $\phi$,
since it uses standard arguments used in DCP.
We note that the sum, scaling, and maximum rule above (\ie, the supremum
rule when $\mathcal A$ is finite) are special cases of DCP-constructed minorants.

\subsection{Memory-based minorants}\label{s-mem}
Suppose we have found minorants $\hat f^i$ of $f$ for iterations $i=1, \ldots, k$.
Then their pointwise maximum
\BEQ\label{e-max-minorant}
\max\{ \hat f^1, \ldots, \hat f^k \}
\EEQ
is also a minorant,
because $\hat f^1, \ldots, \hat f^{k-1}$ are PLBs, 
and $\hat f^k$ is a minorant.
The minorant \eqref{e-max-minorant} is pointwise larger than $\hat f^k$.
We say that the minorant \eqref{e-max-minorant} has \emph{memory}
(of the previously found minorants).
We can also limit the memory to, say, the last $M$ minorants, as
\BEQ\label{e-max-min-mem}
\max \{ \hat f^{k-M}, \ldots, \hat f^k \}.
\EEQ
We refer to this as a limited-memory or finite-memory minorant.

When the minorants in each iteration are affine, we can express
the minorant \eqref{e-max-minorant} as
\BEQ\label{e-pwl-minorant}
\hat{f}^k(x) = \max_{i=1, \ldots, k} \left( f(x^i) + (g^i)^T(x-x^i) \right),
\EEQ
where $g^i \in \partial f(x^i)$.
These minorants are piecewise affine and require only the evaluation of
a subgradient of $f$, and its value, at each point.
When these minorants are used, PMM looks very much like a cutting-plane or
bundle method.

\section{Numerical experiments}\label{s-numexp}
We present two numerical experiments to illustrate the PMM.
Python code for these experiments is available as a
Jupyter notebook at
\begin{quote}
\url{https://github.com/cvxgrp/polyak_minorant}
\end{quote}
The data for both problems was generated with a seeded pseudorandom number generator,
so our numerical experiments can be reproduced exactly.

We used CVXPY \cite{diamond2016cvxpy} 
with the SOCP solver Clarabel \cite{goulart2021clarabel} to compute the
PMM updates.  This incurs some inefficiency, but our goal is only
to illustrate how PMM works with various minorants, with a trade-off of 
per-iteration cost and overall iteration cost.
The numerical examples were run on a
laptop with Ryzen 9 5900HX processor and 32 GB of DDDR4 memory.

\subsection{Second-order cone program}
We consider an instance of the primal-dual cone program,
given in \eqref{e-conepair} and \eqref{e-pdcone-pmm}, with 
\[
\mathcal K = \mathcal K_1 \times \cdots \times \mathcal K_l,
\]
where $\mathcal K_i = \{ (s_i,t_i) \mid \|s\|_2 \leq t \}$
are second-order cones.
These cones are self-dual, \ie, $\mathcal K_i^* = \mathcal K_i$.
In \eqref{e-pdcone-pmm} we list the cone distances separately as 
\[
d_{\mathcal K_i} (x_i)\leq 0, \quad i=1, \ldots, l,
\]
and similarly for the dual cone constraints.
In the form \eqref{e-prob}, we have $m=2l$ inequality constraint functions.
There is an analytical expression for the projection onto a second-order cone,
and from this, we can derive an analytical expression for a subgradient 
of $d_{\mathcal K_i}$.

\paragraph{Data generation.} 
We generate an instance of \eqref{e-conepair}
with $n=500$ variables, $p=200$ equality constraints, and $l=10$ cones each
of dimension $50$.

We generate the data $A$, $b$, and $c$ as follows.
We first generate a vector $z \in \reals^n$ with standard normal entries.
We project $z$ onto $\mathcal K$ to obtain $u$, and we set 
$s=u-z$, which guarantees that $s\in \mathcal K^*=\mathcal K$ 
\cite{rockafellar1981TheTO}[\S 14]. These two vectors satisfy $s^Tu=0$, 
\ie, they are complementary with respect to $\mathcal K$.
Then we generate $A\in \reals^{p \times n}$ and $v\in \reals^p$
with standard normal entries. We set $b=Au$ and $c=s+A^Tv$.
The zero-gap equality constraint $c^Tu=b^Tv$ holds for this data.
For our experiment we only use the data $A$, $b$, and $c$,
and not the primal-dual solution $(u,v,s)$.

\paragraph{Minorant construction and initial point.}
We use a basic subgradient minorant for $d_{\mathcal{K}}$ and
$d_{\mathcal{K}^*}$, and explore different memories $M$.
Our initial point is $x^1=0$.

\begin{figure}
	\centering
	\includegraphics[width=0.95\textwidth]{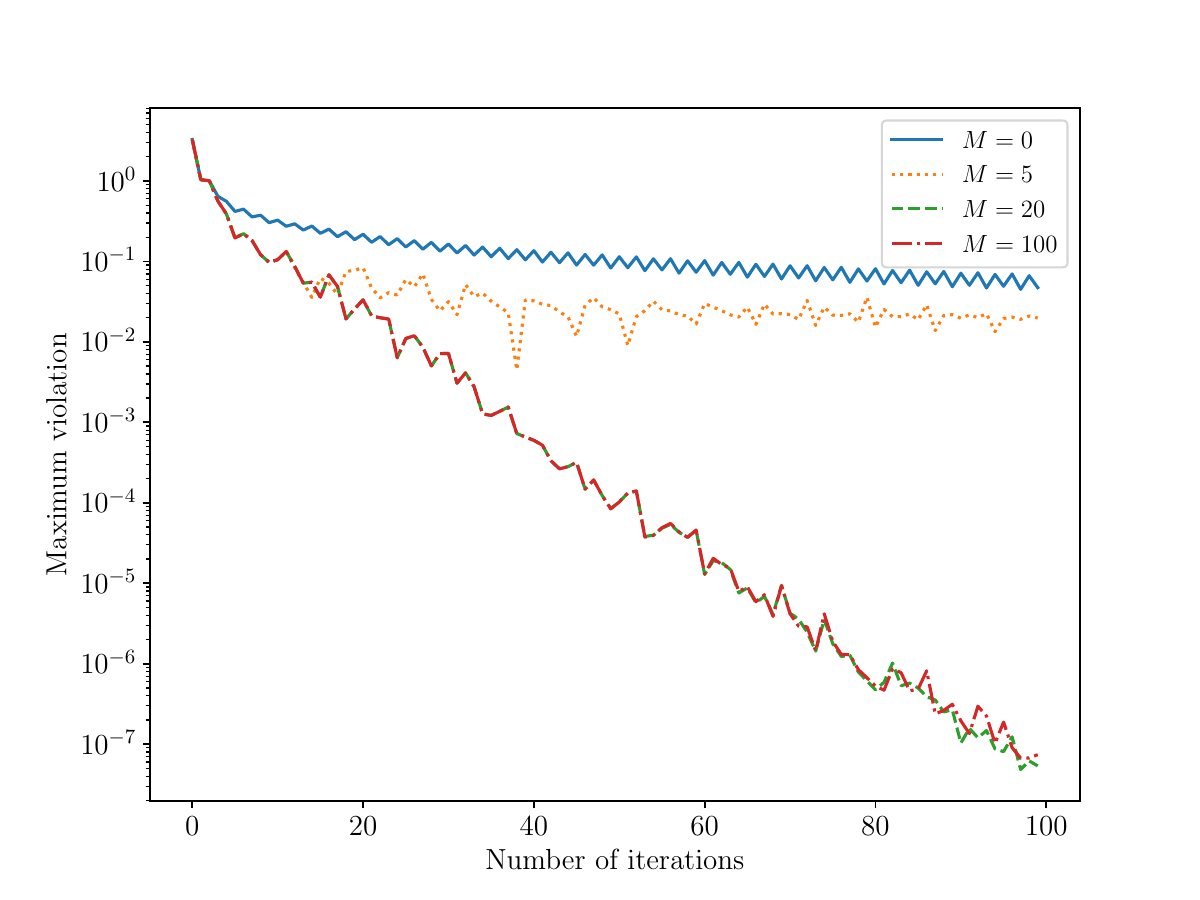}
\caption{Maximum violation versus iterations for primal-dual cone problem.}
	\label{f-socp-iters}
\end{figure}

\begin{figure}
	\centering
	\includegraphics[width=0.95\textwidth]{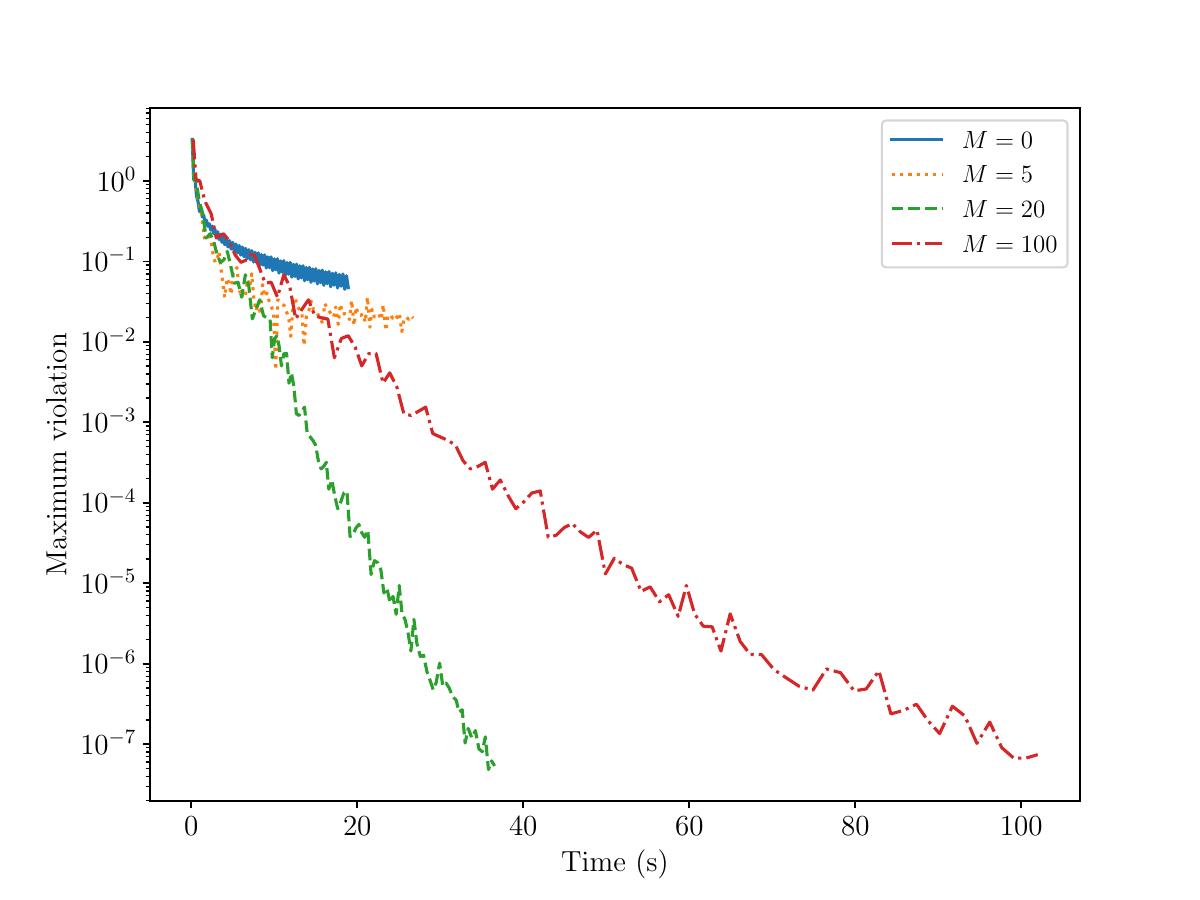}
\caption{Maximum violation versus time for primal-dual cone problem.}
	\label{f-socp-time}
\end{figure}

\paragraph{Results.} 
Figure~\ref{f-socp-iters} shows
the maximum violation $v(x^k)$ versus $k$, the number of iterations,
for memory values $M=0,5,20,100$.
Not surprisingly we get a strong speedup with $M=20$ and $M=100$
compared to smaller memory.
Figure~\ref{f-socp-time} shows
the maximum violation $v(x^k)$ versus the total elapsed time,
which takes into account the varying complexity of the subproblems solved
in each iteration.
Here we see a clear best value of memory $M=20$.

\clearpage 
\subsection{Linear matrix inequality}
Our second example is a linear matrix inequality (LMI).
The goal is to find a matrix $X \in \symm^q$ that satisfies
\[
X \geq I, \quad A_i^T X + XA_i \leq 0, \quad i=1, \ldots, k,
\]
where the inequalities are with respect to the positive semidefinite
cone, and $A_i \in \reals^{q\times q}$, $i=1,\ldots, m$, are given data.
Problems of this form arise in the analysis of control systems; see,
\eg, \cite{boyd1994lmi}.

We express this as the feasibility problem
\BEQ\label{e-lmi-lamb}
\begin{array}{ll}
\mbox{minimize}   & 0 \\
\mbox{subject to} & \lambda_\text{max}(I-X)\leq 0\\
& \lambda_\text{max}(A_i^T X + XA_i)\leq 0, \quad i=1, \ldots, k.
\end{array}
\EEQ
The variable is $X\in \symm^q$, which has dimension $n=q(q+1)/2$.
There are no equality constraints and $m=k+1$ inequality constraints.

\paragraph{Data generation.}
We set $q=20$ and $k=10$, so $n=210$ and $m=11$.
To generate $A_i$ we proceed as follows.
First, we generate 
\[
\tilde A_i = -B_iB_i^T + C_i-C_i^T,
\]
where the entries of the $q \times q$ matrices $B_i$ and $C_i$ are
standard normals.  This means that, with probability one,
\[
\tilde A_i^T + \tilde A_i \leq 0,
\]
\ie, they satisfy the constraints $\tilde A_i^TX+X\tilde A_i\leq 0$ with $X=I$.
Now we generate a $q \times q$ matrix $F$ with
standard normal entries, so $F$ is invertible with probability one,
and form 
\[
A_i = F^{-1}\tilde A_i F.
\]
Then $X=F^TF$ satisfies $A_i^T X+ XA_i \leq 0$,
$i=1, \ldots, k$. Since $X>0$ (with probability one) we can scale it 
to obtain a solution of the LMI \eqref{e-lmi-lamb}.
For our experiment, we only use the data $A_i$, and not the solution $X$.

\paragraph{Minorant construction and initial point.}
Each constraint is the composition of a linear function with $\lambda_\text{max}$.
For $\lambda_\text{max}$ we use the minorant \eqref{e-eigeneigenmin},
with dimension $2$, so the minorants are second-order cone representable,
and the projection can be computed using a SOCP solver. 
We run PMM with memories $M=0,5,20,100$.

\begin{figure}
\centering
\includegraphics[width=0.95\textwidth]{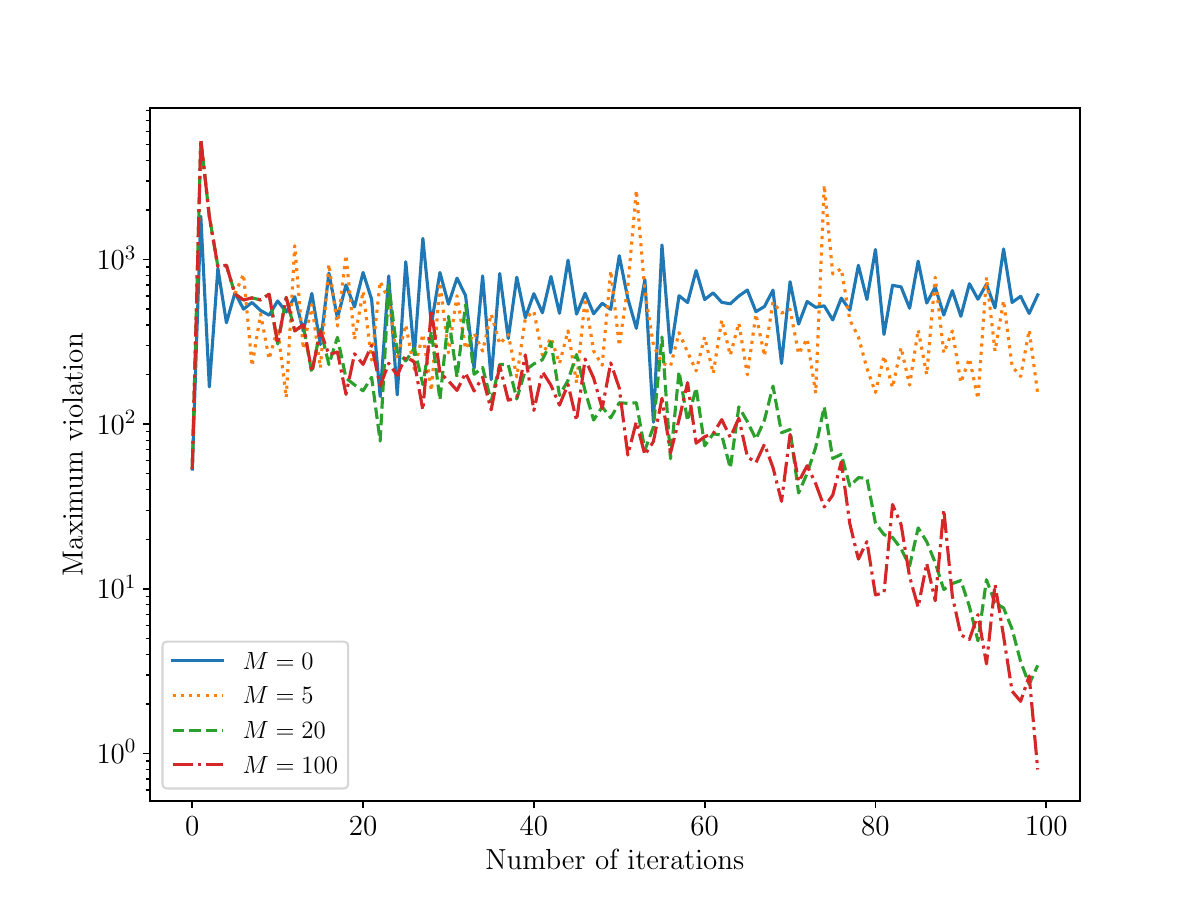}
\caption{Maximum violation versus iterations for LMI problem.}
\label{f-sdp-iters}
\end{figure}

\begin{figure}
\centering
\includegraphics[width=0.95\textwidth]{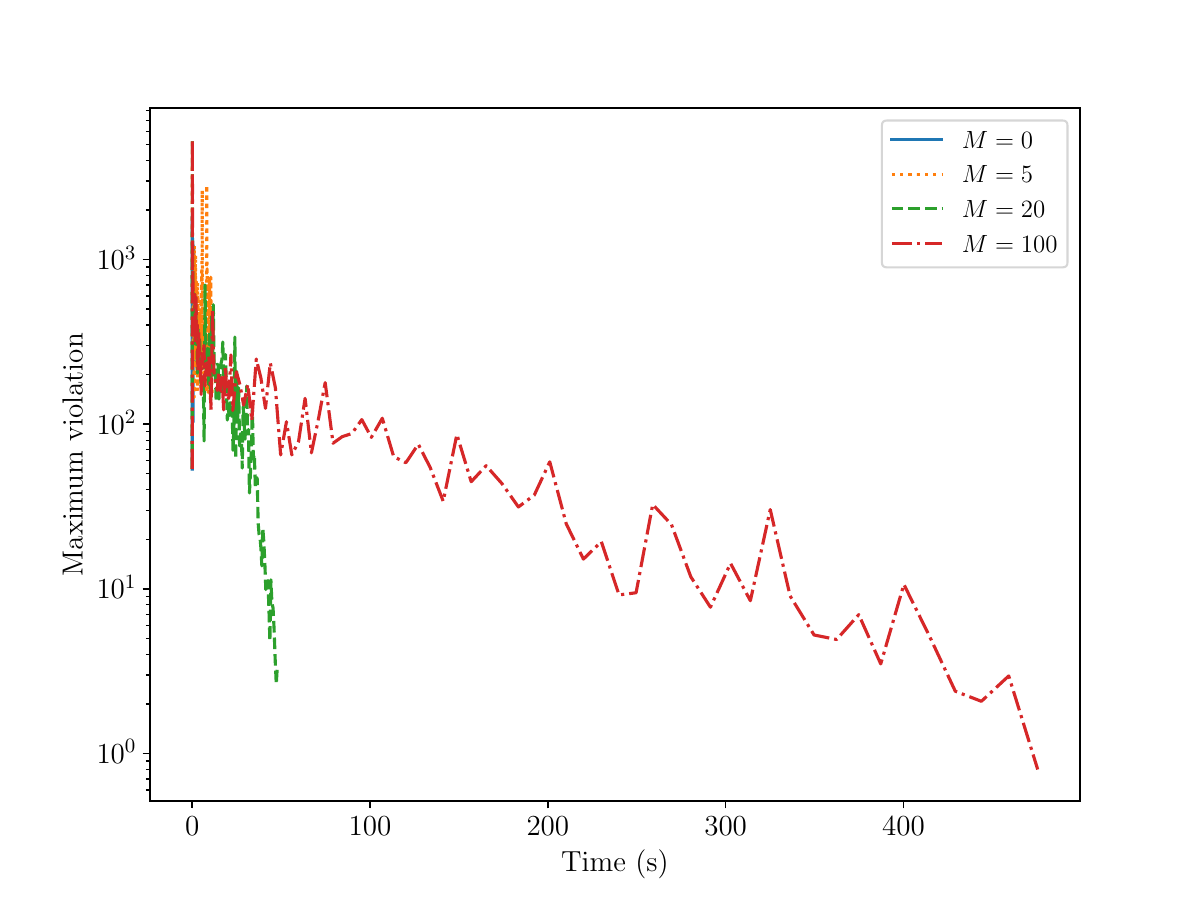}
\caption{Maximum violation versus time for primal-dual cone problem.}
\label{f-sdp-time}
\end{figure}

\paragraph{Results.}
Figure~\ref{f-sdp-iters} shows
the maximum violation $v(x^k)$ versus $k$, the number of iterations,
for memory values $M=0,5,20,100$,
and figure~\ref{f-sdp-time} shows the maximum violation versus elapsed time.
The results are very similar to the previous example, with memory $M=20$ 
giving the fastest (in time) convergence.

\clearpage
\subsection*{Acknowledgments}
This paper builds on notes written around 2010 for the Stanford course
EE364B, \emph{Convex Optimization II}, to which Lieven Vandenberghe, Almir Mutapcic,
Jaehyun Park, Lin Xiao, and Jacob Mattingley contributed.
We thank Tetiana Parshakova, Fangzhao Zhang, Parth Nobel, 
Logan Bell, and Thomas Schmeltzer for useful discussions.
The authors thank an anonymous reviewer for suggesting the 
sharpness convergence result, as well as pointing us to some very 
relevant literature that we had missed in an earlier version of this paper.

Stephen Boyd would like to dedicate this paper to Boris Polyak,
his hero and friend.

\clearpage
\bibliography{refs.bib}

\newcommand{\etalchar}[1]{$^{#1}$}
\begin{thebibliography}{BBTGBT10}

\bibitem[AXKT23]{abdukhakimov2023stochastic}
F.~Abdukhakimov, C.~Xiang, D.~Kamzolov, and M.~Takáč.
\newblock Stochastic gradient descent with preconditioned {P}olyak step-size.
\newblock \url{https://arxiv.org/abs/2310.02093}, 2023.

\bibitem[BBTGBT10]{beck2010comirror}
A.~Beck, A.~Ben-Tal, N.~Guttmann-Beck, and L.~Tetruashvili.
\newblock The comirror algorithm for solving nonsmooth constrained convex
  problems.
\newblock {\em Operations Research Letters}, 38(6):493--498, 2010.

\bibitem[BEGFB94]{boyd1994lmi}
S.~Boyd, L.~El~Ghaoui, E.~Feron, and V.~Balakrishnan.
\newblock {\em Linear Matrix Inequalities in System and Control Theory}.
\newblock Society for Industrial and Applied Mathematics, January 1994.

\bibitem[BTd20]{barre2022complexity}
M.~Barr\'e, A.~Taylor, and A.~d'Aspremont.
\newblock Complexity guarantees for {P}olyak steps with momentum.
\newblock In {\em Proceedings of Thirty Third Conference on Learning Theory},
  volume 125, pages 452--478. PMLR, 09--12 Jul 2020.

\bibitem[BV04]{boyd2004convex}
S.~Boyd and L.~Vandenberghe.
\newblock {\em Convex Optimization}.
\newblock Cambridge University Press, 2004.

\bibitem[BV07]{boyd2007localization}
S.~Boyd and L.~Vandenberghe.
\newblock Localization and cutting-plane methods.
\newblock Lecture notes for EE364b, Stanford University, 2007.

\bibitem[BZK21]{berrada2021comment}
L.~Berrada, A.~Zisserman, and M.~Kumar.
\newblock Comment on stochastic {P}olyak step-size: {P}erformance of {ALI-G}.
\newblock \url{https://arxiv.org/abs/2105.10011}, 2021.

\bibitem[CG59]{cheney1959newton}
E.~Cheney and A.~Goldstein.
\newblock Newton's method for convex programming and {T}chebycheff
  approximation.
\newblock {\em Numerische Mathematik}, 1:253--268, 1959.

\bibitem[CL12]{cheng2012active}
W.~Cheng and D.~Li.
\newblock An active set modified
  {P}olak{\textendash}{R}ibi{\'{e}}re{\textendash}{P}olyak method for
  large-scale nonlinear bound constrained optimization.
\newblock {\em Journal of Optimization Theory and Applications},
  155(3):1084--1094, June 2012.

\bibitem[DB16]{diamond2016cvxpy}
S.~Diamond and S.~Boyd.
\newblock {CVXPY}: {A} {P}ython-embedded modeling language for convex
  optimization.
\newblock {\em Journal of Machine Learning Research}, 17(83):1--5, 2016.

\bibitem[DDMP18]{davis2018subgradient}
D.~Davis, D.~Drusvyatskiy, K.~MacPhee, and C.~Paquette.
\newblock Subgradient methods for sharp weakly convex functions.
\newblock {\em Journal of Optimization Theory and Applications},
  179(3):962–982, September 2018.

\bibitem[Dru20]{drusvyatskiy2020convex}
D.~Drusvyatskiy.
\newblock Convex analysis and nonsmooth optimization.
\newblock
  \url{https://sites.math.washington.edu/~ddrusv/crs/Math_516_2020/bookwithindex.pdf},
  2020.

\bibitem[Fra20]{frangioni2020standard}
A.~Frangioni.
\newblock Standard bundle methods: {U}ntrusted models and duality.
\newblock In {\em Numerical Nonsmooth Optimization}, pages 61--116. Springer,
  2020.

\bibitem[GBGP22]{gower2022cutting}
R.~Gower, M.~Blondel, N.~Gazagnadou, and F.~Pedregosa.
\newblock Cutting some slack for {SGD} with adaptive {P}olyak stepsizes.
\newblock \url{https://arxiv.org/abs/2202.12328}, 2022.

\bibitem[GBY06]{grant2006disciplined}
M.~Grant, S.~Boyd, and Y.~Ye.
\newblock Disciplined convex programming.
\newblock In {\em Global Optimization}, pages 155--210. Springer, 2006.

\bibitem[GC21]{goulart2021clarabel}
P.~Goulart and Y.~Chen.
\newblock Clarabel: {A} library for optimization and control, 2021.
\newblock URL: {https://oxfordcontrol.github.io/ClarabelDocs/stable/}.

\bibitem[Gof77]{goffin1977rates}
J.~Goffin.
\newblock On convergence rates of subgradient optimization methods.
\newblock {\em Mathematical Programming}, 13(1):329–347, December 1977.

\bibitem[GTD22]{goujaud2022quadratic}
B.~Goujaud, A.~Taylor, and A.~Dieuleveut.
\newblock Quadratic minimization: {F}rom conjugate gradients to an adaptive
  heavy-ball method with {P}olyak step-sizes.
\newblock In {\em OPT 2022: Optimization for Machine Learning (NeurIPS 2022
  Workshop)}, 2022.

\bibitem[HK22]{hazan2022revisiting}
E.~Hazan and S.~Kakade.
\newblock Revisiting the {P}olyak step size.
\newblock \url{https://arxiv.org/abs/1905.00313}, 2022.

\bibitem[HUL96]{hiriart1996convex}
J.~Hiriart-Urruty and C.~Lemar{\'e}chal.
\newblock {\em Convex Analysis and Minimization Algorithms {II}: {A}dvanced
  Theory and Bundle Methods}.
\newblock Grundlehren der mathematischen Wissenschaften. Springer Berlin
  Heidelberg, 1996.

\bibitem[Jon98]{jongbloed1998the}
G.~Jongbloed.
\newblock The iterative convex minorant algorithm for nonparametric estimation.
\newblock {\em Journal of Computational and Graphical Statistics}, 7(3):310,
  September 1998.

\bibitem[Kel60]{kelley1960cutting}
J.~Kelley.
\newblock The cutting-plane method for solving convex programs.
\newblock {\em Journal of the Society for Industrial and Applied Mathematics},
  8(4):703--712, 1960.

\bibitem[Kiw90]{kiwiel1990proximity}
K.~Kiwiel.
\newblock Proximity control in bundle methods for convex nondifferentiable
  minimization.
\newblock {\em Mathematical Programming}, 46(1-3):105--122, 1990.

\bibitem[Kow09]{kowalewski1909power}
G.~Kowalewski.
\newblock {\em Einf{\"u}hrung in die determinantentheorie einschliesslich der
  unendlichen und der Fredholmschen determinanten}.
\newblock Veit \& comp., 1909.

\bibitem[LMY18]{lin2018level}
Q.~Lin, R.~Ma, and T.~Yang.
\newblock Level-set methods for finite-sum constrained convex optimization.
\newblock In {\em Proceedings of the 35th International Conference on Machine
  Learning}, volume~80, pages 3112--3121. PMLR, 10--15 Jul 2018.

\bibitem[LNN95]{lemarechal1995new}
C.~Lemar{\'e}chal, A.~Nemirovskii, and Y.~Nesterov.
\newblock New variants of bundle methods.
\newblock {\em Mathematical Programming}, 69(1):111--147, 1995.

\bibitem[LST{\etalchar{+}}23]{li2023sp2}
S.~Li, W.~Swartworth, M.~Takáč, D.~Needell, and R.~Gower.
\newblock {SP}2 : {A} second order stochastic {P}olyak method.
\newblock In {\em The Eleventh International Conference on Learning
  Representations}, 2023.

\bibitem[LVHLLJ21]{loizou2021stochastic}
N.~Loizou, S.~Vaswani, I.~Hadj~Laradji, and S.~Lacoste-Julien.
\newblock Stochastic {P}olyak step-size for {SGD}: {A}n adaptive learning rate
  for fast convergence, 13--15 Apr 2021.

\bibitem[LZ20]{lan2020algorithms}
G.~Lan and Z.~Zhou.
\newblock Algorithms for stochastic optimization with function or expectation
  constraints.
\newblock {\em Computational Optimization and Applications}, 76(2):461–498,
  February 2020.

\bibitem[McL78]{mclinden1978affine}
L.~McLinden.
\newblock Affine minorants minimizing the sum of convex functions.
\newblock {\em Journal of Optimization Theory and Applications},
  24(4):569--583, April 1978.

\bibitem[MHB75]{marsten1975boxstep}
R.~Marsten, W.~Hogan, and J.~Blankenship.
\newblock The boxstep method for large-scale optimization.
\newblock {\em Operations Research}, 23(3):389--405, 1975.

\bibitem[Nem96]{nemirovski1996interior}
A.~Nemirovski.
\newblock Interior point polynomial time methods in convex programming.
\newblock \url{https://www2.isye.gatech.edu/~nemirovs/Lect_IPM.pdf}, 1996.

\bibitem[Nes18]{nesterov2018lectures}
Y.~Nesterov.
\newblock {\em Lectures on Convex Optimization}.
\newblock Springer, 2018.

\bibitem[NN94]{nesterov1994interior}
Y.~Nesterov and A.~Nemirovskii.
\newblock {\em Interior-Point Polynomial Algorithms in Convex Programming}.
\newblock Society for Industrial and Applied Mathematics, January 1994.

\bibitem[OLJL22]{orvieto2022dynamics}
A.~Orvieto, S.~Lacoste-Julien, and N.~Loizou.
\newblock Dynamics of {SGD} with stochastic polyak stepsizes: Truly adaptive
  variants and convergence to exact solution.
\newblock In {\em Advances in Neural Information Processing Systems}, 2022.

\bibitem[PB14]{parikh2014proximal}
N.~Parikh and S.~Boyd.
\newblock Proximal algorithms.
\newblock {\em Foundations and Trends in Optimization}, 1(3):127--239, 2014.

\bibitem[PO21]{prazeres2021stochastic}
M.~Prazeres and A.~Oberman.
\newblock Stochastic gradient descent with polyak's learning rate.
\newblock {\em Journal of Scientific Computing}, 89(1), September 2021.

\bibitem[Pol63]{polyak1963gradient}
B.~Polyak.
\newblock Gradient methods for the minimisation of functionals.
\newblock {\em {USSR} Computational Mathematics and Mathematical Physics},
  3(4):864--878, January 1963.

\bibitem[Pol67]{polyak1967general}
B.~Polyak.
\newblock A general method of solving extremum problems.
\newblock {\em Sov. Math., Dokl.}, 8:593--597, 1967.

\bibitem[Pol69]{polyak1969minimization}
B.~Polyak.
\newblock Minimization of unsmooth functionals.
\newblock {\em USSR Computational Mathematics and Mathematical Physics},
  9(3):14–29, January 1969.

\bibitem[Pol87]{polyak1987intro}
B.~Polyak.
\newblock {\em Introduction to optimization}.
\newblock Optimization Software, Inc., 1987.

\bibitem[PZB23]{parshakova2022implementation}
T.~Parshakova, F.~Zhang, and S.~Boyd.
\newblock Implementation of an oracle-structured bundle method for distributed
  optimization.
\newblock {\em Optimization and Engineering}, 2023.

\bibitem[Roc81]{rockafellar1981TheTO}
R.~Rockafellar.
\newblock {\em The Theory of Subgradients and its Applications to Problems of
  Optimization}.
\newblock Heldermann Verlag, 1981.

\bibitem[SBG{\etalchar{+}}20]{stellato2020osqp}
B.~Stellato, G.~Banjac, P.~Goulart, A.~Bemporad, and S.~Boyd.
\newblock {OSQP}: an operator splitting solver for quadratic programs.
\newblock {\em Mathematical Programming Computation}, 12(4):637--672, 2020.

\bibitem[Sho73]{shor1973convergence}
N.~Shor.
\newblock Convergence rate of the gradient descent method with dilatation of
  the space.
\newblock {\em Cybernetics}, 6(2):102–108, 1973.

\bibitem[Sho12]{shor2012minimization}
N.~Shor.
\newblock {\em Minimization Methods for Non-differentiable Functions},
  volume~3.
\newblock Springer Science \& Business Media, 2012.

\bibitem[SNW12]{sra2012optimization}
S.~Sra, S.~Nowozin, and S.~Wright.
\newblock {\em Optimization for Machine Learning}.
\newblock MIT Press, 2012.

\bibitem[vAFdO16]{van2016inexact}
W.~van Ackooij, A.~Frangioni, and W.~de~Oliveira.
\newblock Inexact stabilized {B}enders’ decomposition approaches with
  application to chance-constrained problems with finite support.
\newblock {\em Computational Optimization and Applications}, 65:637--669, 2016.

\bibitem[WJZ23]{wang2023generalized}
X.~Wang, M.~Johansson, and T.~Zhang.
\newblock Generalized {P}olyak step size for first order optimization with
  momentum.
\newblock \url{https://arxiv.org/abs/2305.12939}, 2023.

\bibitem[YCL22]{you2022minimizing}
J.~You, H.~Cheng, and Y.~Li.
\newblock Minimizing quantum {R}ényi divergences via mirror descent with
  {P}olyak step size.
\newblock In {\em 2022 IEEE International Symposium on Information Theory
  (ISIT)}, pages 252--257, 2022.

\bibitem[YL22]{you2022two}
J.~You and Y.~Li.
\newblock Two {P}olyak-type step sizes for mirror descent.
\newblock \url{https://arxiv.org/abs/2210.01532}, 2022.

\end{thebibliography}

\end{document}